\documentclass[10pt,
]{amsart}
\usepackage{
amsmath, amsthm, amssymb, graphicx, commath, showkeys, pgfplots, pst-solides3d, inputenc, float, amssymb, graphicx, verbatim, tikz-cd
}
\newcommand{\no}[1]{#1}
\renewcommand{\no}[1]{}  \newcommand{\upDelta}{\Delta} 

\no{\usepackage{times}\usepackage[
subscriptcorrection, slantedGreek, 
nofontinfo]{mtpro}
\renewcommand{\Delta}{\upDelta}
}
\numberwithin{equation}{section}% 

\newtheorem{theorem}{Theorem}[section]
\newtheorem{proposition}{Proposition}[section]
\newtheorem{lemma}{Lemma}[section]
\newtheorem{definition}{Definition}[section]
\newtheorem{corollary}{Corollary}[section]

\theoremstyle{definition}
\newtheorem{remark}{Remark}[section]

\DeclareMathOperator{\supp}{supp}

\DeclareMathOperator{\WF}{WF}
\DeclareMathOperator{\WFA}{WF_A}

\DeclareMathOperator{\n}{neigh}

\DeclareMathOperator{\sgn}{sgn}

\newcommand{\R}{{\bf R}}

\renewcommand{\r}[1]{(\ref{#1})}

\newcommand{\be}[1]{\begin{equation}\label{#1}}
\newcommand{\ee}{\end{equation}}
\renewcommand{\d}{\mathrm{d}}

\renewcommand{\i}{\mathrm{i}}

\title[Invertibility and Stability for A Generic Class of Radon Transform]{Invertibility and Stability for A Generic Class of Radon Transforms with Application to Dynamic Operators}

\author[S.Rabienia]{Siamak RabieniaHaratbar}
\address{Department of Mathematics, Purdue University, West Lafayette, IN 47907}
\thanks{ Partly supported by NSF Grant DMS ~1600327}

\begin{document}
\begin{abstract}
	Let $X$ be an open subset of $\R^2$. We study the dynamic operator, $\mathcal{A}$, integrating over a family of level curves in $X$ when the object changes between the measurement. We use analytic microlocal analysis to determine which singularities can be recovered by the data-set. Our results show that not all singularities can be recovered, as the object moves with a speed lower than the X-ray source. We establish stability estimates and prove that the injectivity and stability are of a generic set if the dynamic operator satisfies the visibility, no conjugate points, and local Bolker conditions. We also show this results can be implemented to Fan beam geometry.  
\end{abstract}

\maketitle
	
\section{Introduction}
	Tomography of moving objects has been attracting a growing interest recently, due to its wide range of applications in medical imaging, for example, X-ray of the heart or the lungs. Data acquisition and reconstruction of the object which changes its shape during the measurement is one of the challenges in computed tomography and dynamic inverse problems. The major difficulty in the reconstruction of images from the measurement sets is the fact that object changes between measurements but does not move fast enough compared to the speed of X-rays. This means that some singularities of the object might not be detectable even if the source fully rotates around the object. The application of known reconstruction methods (based on the inversion of the Radon transform) usually results in many motion artifacts within the reconstructed images if the motion is not taken into account. One extreme example will be the case when the object (or some small part of it) rotates with the same rate as the scanner. This leads to integration over the same family of rays (see also \cite{25}), and therefore, one cannot locally recover all the singularities.     
	
	Analytic techniques for reconstruction of dynamic objects, known as motion compensation, have been used widely for different types of motion, like affine deformation, see e.g. \cite{5, 6, 12, 13, 18, 20, 21, 26}. In the case of non-affine deformations, there is no inversion formula. Iterative reconstructions, however, do exist in order to detect singularities by approximation of inversion formulas for the parallel and fan beam geometries \cite{19}, as well as cone beam geometry \cite{22}. In a recent work, Hahn and Quinto \cite{11} studied the dynamic operator 
\be{1.1}
	\mathcal{A} f(s,t)=\int_{z\cdot \omega(t)=s} \mu(t,z) f(\psi_t (z)) \d S_z, \quad \text{  $\omega(t)=(\cos t, \sin t)$,}
\ee
	with a smooth motion where the limited data case has been analyzed, and characterization of visible and added singularities have been investigated.
	
	Our work in this paper is motivated by these dynamic measurements. We first show this dynamic problem can be reduced to an integral geometry problem integrating over level curves. By an appropriate change of variable (see section 2), $\mathcal{A}$ can be written as 
$$
	\mathcal{A} f(s,t)=\int_{\psi_t ^{-1}(x)\cdot \omega(t)=s} \hat{\mu}(t,x) f(x) \d S.
$$
	Therefore, we study the following general operator:
$$
	\mathcal{A} f(s,t)=\int_{\phi(t,x)=s} {\mu}(t,x) f(x) \d S_{s,t},
$$
	which allows us to study the original dynamic problem with a more general set of curves (see also \cite{7},) and then transfer the result to a dynamic operator $\mathcal{A}$ given by \r{1.1}.
			
	The dynamic operator $\mathcal{A}$ formulated as above falls into the general microlocal framework studied by Beylkin \cite{1} (see also \cite{13}) which goes back to Guillemin and Sternberg \cite{9, 10} who studied the integral geometry problems with a more general platform from the microlocal point of view. See also \cite{7}, where a weighted integral transform has been studied on a compact manifold with a boundary over a general set of curves (a smooth family of curves passing through every point in every direction).     
	
	The main novelty of our work, compared to previous works which are concentrated on the microlocal invertibility, is that for the dynamic problem, under some natural microlocal conditions, the actual uniqueness and stability results have been established. In fact, our imposed natural microlocal conditions guarantee that one can recover each singularity, and a functional analysis argument leads to stability results. We show that under these conditions, the dynamic operator is stably invertible in a neighborhood of pairs $(\phi,\mu)$ in a generic set, and in particular, it is injective and stable for slow enough motion (which is not required to be a periodic motion model). This is the similar kind of stability result which has been studied in \cite{17} for the generalized Radon transform and in \cite{7} which coincide when the dimension is two. The data is cut (restricted) in a way to have the normal operator related to the localized dynamic operator $\mathcal{A}$ as a pseudodifferential operator ($\Psi$DO) near each singularity. We do not analyze the case where these conditions are not satisfied globally, but our analysis (see also \cite{11}) shows that one can still recover the visible singularities in a stable way, and periodicity or non-periodicity plays no role in the reconstruction process. We also show that, due to the generality of our approach, our results can be implemented to other geometries, for instance, fan beam geometry. 
	
	This paper is organized as follow: Section one is an introduction. In section two, we state the definitions of $\it{Visibility}$, $\it{Local}$ $\it{Bolker}$ $\it{Condition}$, $\it{Semi}$-$\it{Global}$ $\it{Bolker}$ $\it{Condition}$, and our main result. Some preliminary results have been stated in section three. Section four is devoted to analytic microlocal analysis approach which is used to show that the operator $\mathcal{A}$ is a Fourier Integral Operator (FIO). Then the canonical relation $\mathcal{C}$ is computed and it is shown that it is a four-dimensional non-degenerated conic submanifold of the conormal bundle. In section five, it is shown that a certain localized version of the normal operator $\mathcal{N}=\mathcal{A}^*\mathcal{A}$ is an elliptic pseudodifferential operator ($\Psi$DO) under the visibility, and the local and semi-global Bolker conditions. In section six, we study the operators $\mathcal{A}$ and $\mathcal{N}$ globally, and show that uniqueness and stability (injectivity) are of a generic set with the corresponding topology. In the last section, we implement our results for the initial dynamic problem of scanning a moving object while changing its shape. We also show that our results can be applied to fan beam geometry by an appropriate choice of phase function $\phi.$ 
		   
\section{Main Results}
	In this section, we first introduce the dynamic operator and then reduce it to an integral geometry problem integrating over level curves. After some necessary propositions, we state our main results.
\begin{definition}
	Let $X$ be a fixed open set in $\R^2$ and $Y$ be the open sets of lines determined by $(s,t)$ in $\R^2$. For $\mathcal{A}: C_0 ^\infty (X) \rightarrow C ^\infty (Y)$, the operator of the dynamic inverse problem is defined by
$$
	\mathcal{A} f(s,t)=\int_{x\cdot \omega(t)=s} \mu(t,x) f(\psi_t (x)) \d S_x,
$$
	where $\omega(t)=(\cos t,\sin t)$ and the function $\mu$ is a non-vanishing smooth weight changing with respect to the variable $t$ and the position $x$.
\end{definition}
	Here $\psi_t$ is a diffeomorphism in $\R^2$, which is identity outside $X$, smoothly depending on the variable $t$, and $\d S_x$ is the euclidean measure restricted to the lines parametrized by $\{s=x\cdot\omega(t)\}$. Notice that each point (position) $x\in X$, lies on the lines in $Y$ parametrized by $(s,t)$.

	The operator $\mathcal{A}$ can be written in the following format:
$$
	\mathcal{A} f(s,t)=\iint_{\R^2} \mu(t,x) f(\psi_t (x))\delta(s-x\cdot\omega(t)) \d x.
$$
	Since $\psi_t$ is a diffeomorphism, by performing a change of variable $z=\psi_t (x)$, we get $x=\psi_t ^{-1}(z)$ and therefore, we have
$$
	\mathcal{A} f(s,t)=\iint_{\R^2} J(t,z) \mu(t,\psi_t ^{-1}(z)) f(z)\ \delta(s-\psi_t ^{-1}(z)\cdot \omega(t)) \ dz. 
$$
	From now on, we do most of our analysis on the following general operator:
\be{2.1}
	\mathcal{A} f(s,t)=\int_{\phi(t,x)=s} {\mu}(t,x) f(x) \d S_{s,t},
\ee
	where $\mu$ is a new positive and real analytic weight and the map  
$$
	x=(x^1, x^2)\longrightarrow \phi(t,x),
$$
	with analytic function $\phi$, is real-valued. Here $\d S_{s,t}$ is the Euclidean measure of the level curves of function $\phi$, defined as
$$
	H(s,t)=\lbrace x\in X : s=\phi(t,x) \rbrace, \quad  s\in \R, \  t\in \R.
$$ 
	We, first, need to show for any time $t$ and point $x$, there exists a curve passing through the point $x$ with direction $\omega(t)$. 
\begin{proposition}
	Let $H(s,t)$ be the level curves of $\phi$. Then, locally near $(s_0,t_0)$ and near a fixed $x_0\in H(s,t)$ the followings are equivalent.
	\\i) The map from the variable $t$ to the unit normal vector $\nu$ of the level curves $H(s,t)$: 
\be{2.2}
	t\longrightarrow \nu(t,x)=\frac{\partial_x \phi(t,x)}{|\partial_x \phi(t,x)|}, \quad \partial_x \phi(t,x)\not = 0,
\ee	
	is a local diffeomorphism, where $\partial_x=(\partial_{x^1},\partial_{x^2})$.
	\\ii) The {\bfseries{Local Bolker Condition}}: 
\be{2.3}
	h(t,x)=\det \left( {\begin{array}{cc}
	\frac{\partial \phi}{ \partial x^j},  \frac{\partial^2 \phi}{\partial t \partial x^j}       \end{array} } \right)_{\big|_{(t,x)=(t_0,x_0)}}\not = 0,
\ee
	holds locally near $(s_0,t_0)$ and near $x_0$.
\end{proposition}
\begin{remark}
	i) The proof of Proposition 2.1 is postponed to the next section. In our setting, the equation \r{2.3} is the generalization of what it is known as a Bolker condition in [Theorem~ 14 (2), \cite{11}].\vspace{1.5mm}\\	
	ii) One can always rotate the unit normal vector $\nu$ by $\frac{\pi}{2}$ (at a fixed point $x$ on the curve) to get the tangent vector at that fixed point. Now the first part in Proposition 2.1 implies that the map from the variable $t$ to the tangent vector at point $x$ on the level curve $H(s,t)$, is also a local diffeomorphism.\vspace{1.5mm}\\	
	iii) We work locally near $(s_0,t_0)$ and a fixed $x_0$ on the level curve. Let $l_0$ denote the unit tangent (normal) vector at $x_0$. By the first part, for any unit tangent vector $l$ in some small neighborhood of $l_0$ ($l$ is some perturbation of $l_0$), the map from the variable $t$ to the unit tangent vector at a fixed point $x$ is a local diffeomorphism. Now the Implicit Function Theorem implies that for any given $t$, there exists a curve passing through the fixed point $x$ with a tangent vector $l$. This indeed is what to expect if we want the level curves to behave like the geodesic curves.\vspace{1.5mm}\\	
	iv) The local Bolker condition requires that when the object moves in time, the curve changes its direction. A counterexample when the local Bolker condition does not hold is the case where an object and the scanner move with the same rate. In this situation, the object can be considered stationary where it is being scanned with stationary parallel rays. The above proposition guarantees that locally and microlocally this situation will not happen and the parameter $t$ changes the angle if we keep the object stationary. (i.e the movement is not going to be synchronized with the scanner)\vspace{1.5mm}\\	
	v) Proposition 3.1 in the next section, shows that one can connect the local Bolker condition to Fourier Integral Operator (FIO) theory by extending the function $\phi$ to a homogeneous function of order one (see \cite{1}), and therefore one can use the condition \r{2.3} for the analysis. 
\end{remark}

For main results, we first state the following definitions.
	
\begin{definition}
	The function $\phi$ satisfies the {\bfseries{Visibility condition}} at $(x,\xi) \in T^*X \setminus 0$ if there exists a pair $(s,t)$ with property $\phi(t,x)=s,$ such that $\partial_{x}\phi(t,x) \parallel \xi.$ Here $T^*X$ is the cotangent bundle of $X$.
\end{definition}
	The visibility condition requires that at a point $x$ and co-direction $\xi$, locally, there exists a curve passing through $x$ which is conormal to $\xi$. As we pointed out in Remark 2.1, this property is a natural property of level curves as are expected to behave like geodesic curves. It also means that each singularity can be probed locally.
\begin{definition}
	Let $(x_0,(s_0,t_0))\in X\times Y$ be a fixed point with property $s_0=\phi(t_0,x_0)$. The function $\phi$ satisfies the {\bfseries{Semi-Global Bolker Condition}} at $(x_0,(s_0,t_0))$ if there exists a neighborhood of $(x_0,(s_0,t_0))$, $V$ and $U$, such that for any $(x,(s,t))\in V\times U$ and $y \in X$
\be{2.4}
	\left\{\begin{array}{ll}
		\phi (t, x)=\phi (t, y)=s\\
		\partial_t\phi (t, x)=\partial_t\phi (t, y)\\
	\end{array}
	\right. \ \Longrightarrow \ x=y. %\quad  \quad \quad \text{for $\ y \in X$.}
\ee
\end{definition}
	The first equation in \r{2.4} implies that at instance $t$, both points $x$ and $y$ belong to the same level curve $\phi$. The second equation implies that a perturbation in the variable $t$, cannot distinguish between these two points as they both belong to the same perturbed level curve. Note that, if the level curves $\phi$ are geodesics, it is required that the point $x$ (close to a fixed point $x_0$) has no conjugate points along the curve passing through it with conormal $\xi$. This is indeed a semi-global condition, as $x$ only varies in the open set $V$, but $y$ can be anywhere along the level curve $\phi$, not necessary close to $x$.  
		
	We now are ready to state our main result for the operator $\mathcal{A}$ given by \r{2.1}.
\begin{theorem}
	Consider the operator $\mathcal{A}$ with a nowhere vanishing smooth weight $\mu$. Let $\Sigma$ be a set of all possible pairs $(\phi,\mu)$ which are smooth in some $C^k$-topology with $k$ an arbitrary large natural number. Assume that for any $(x,\xi)\in T^*X\setminus0$, (i) the visibility condition holds and (ii) the local and semi-global Bolker conditions are satisfied for some $(s,t)$ given by the visibility condition. 
	\vspace{2mm}\\
	Then within $\Sigma$, there exists a dense and open (generic) set $\Lambda$ of pairs of $(\phi,\mu)$ such that locally near any pair in $\Lambda$, the uniqueness results and therefore stability (injectivity) estimates given by Proposition 6.2 hold.
\end{theorem}

	To formulate above result for the dynamic operator $\mathcal{A}$ given by \r{1.1}, we first state the visibility, and the local and semi-global Bolker conditions for $\mathcal{A}$.
	\vspace{2mm}
		
	{\it{Visibility.}} This condition implies that for $(x,\xi)\in T^*X \setminus 0,$ the map 
\be{2.5}
	t \rightarrow \frac{\xi}{|\xi|} \in S^1 
\ee
	is locally surjective. Here the point $(s,t)$ lies on the level curve $s=\psi_t ^{-1}(x)\cdot\omega(t).$ 
	\vspace{2mm}
	
	{\it{Local Bolker Condition.}} 
	This condition (see Proposition 4.1) implies that 
\be{2.6}
	h(t,x) = \det \left( {\begin{array}{cc} 
		\frac{\partial \psi_t ^{-1}(x)\cdot\omega(t)}{ \partial x^j},  \frac{\partial^2 \psi_t ^{-1}(x)\cdot\omega(t)}{\partial t \partial x^j} \end{array} } \right)\not=0. 
\ee
	\vspace{2mm}
	{\it{Semi-global Bolker condition (No conjugate points condition).}} By condition \r{2.4}, semi-global Bolker condition holds if the map
\be{2.7}
	x \rightarrow \left( \psi_t ^{-1}(x)\cdot\omega(t),  \partial_t (\psi_t ^{-1}(x)\cdot\omega(t)) \right)
\ee
	is one-to-one.
	
	Now for the dynamic forward operator $\mathcal{A}$ given by \r{1.1}, we have the following result:	
	
\begin{theorem}
	Consider the dynamic operator $\mathcal{A}$ with a nowhere vanishing smooth weight $\mu$. Let $\Sigma$ be a set of all possible pairs $(\psi,\mu)$ which are smooth in some $C^k$-topology with $k$ an arbitrary large natural number. Assume that for any $(x,\xi)\in T^*X\setminus 0$, (i) the visibility condition \r{2.5} holds and (ii) the local and semi-global Bolker conditions given by \r{2.6} and \r{2.7} are satisfied for some $(s,t)$ given by the visibility condition.\vspace{2mm}\\		
	Then within $\Sigma$, there exists a dense and open (generic) set $\Lambda$ of pairs of $(\psi,\mu)$ such that locally near any pair in $\Lambda$, the uniqueness results and therefore stability (injectivity) estimates hold.
\end{theorem}
	
\begin{corollary}
	In particular, for a small perturbation of $\phi(t,x)=x\cdot\omega(t)$ where there is no motion or the motion is small enough $(\mu\approx 1)$, we have the actual injectivity and invertibility as the set of pairs of $(\phi,\mu)$ is included in $\Lambda.$ 
\end{corollary}
\begin{remark}
	The Corollary 2.1 follows from the fact that the stationary Radon transform is analytic and for a small perturbation of phase function, the invertibility and injectivity still hold. 
\end{remark}
	
\section{Preliminary Results}
	In this section, we first prove Proposition 2.1 and then connect the local Bolker condition \r{2.3} to Fourier Integral Operator theory. At the end, we state some definitions which will be used in the following sections.
\begin{definition}
	A set $\Sigma$ is conic, if $\xi\in\Sigma$ then $r\xi\in\Sigma$ for all $r > 0$.
\end{definition}	
\begin{proof}[\textbf{Proof of Proposition 2.1}]
	i) $\rightarrow$ ii) Fix $(t_0,x_0)$ and let $\phi(t_0,x_0)=s_0$. We work on some neighborhood of $(s_0,t_0)$ and $x_0$. Since $\partial_x \phi(t,x)\not = 0$, the map \r{2.2} is well-defined and there exists a tangent at a fixed time $t$ when $x$ varies. The map \r{2.2} is a local diffeomorphism, therefore $\partial_t \nu(t,x)\not=0$ and its inverse exists with non-zero derivative in a conic neighborhood. 
	
	Assume now that $h(t,x)=0$. Then there exists a non-zero constant $c$ such that  
\be{3.1}
	\partial_t \partial_x \phi(t,x)=c\partial_x \phi(t,x). 
\ee
	Plugging \r{3.1} into $\partial_t \nu(t,x)$:  
$$
	\partial_t \nu(t,x)= \frac{\partial_t \partial_x \phi(t,x)}{|\partial_x \phi(t,x)|} - \partial_x \phi(t,x)\frac{\partial_x \phi(t,x)\cdot \partial_t\partial_x \phi(t,x)}{|\partial_x \phi(t,x)|^3}
$$
	we get $\partial_t \nu(t,x)=0$, which is a contradiction. Therefore  
$$
	h(t,x)\not = 0.
$$  
	ii)$\rightarrow$ i) Assume that \r{2.3} is true. This in particular implies that $\partial_x \phi(t,x)$ and $\partial_t \partial_x \phi(t,x)$ are non-zero and linearly independent. For any $t$, let $\nu(t,x)=\frac{\partial_x \phi(t,x)}{|\partial_x \phi(t,x)|}$ denotes the unit normal at a fixed point $x$ on the curve. To show the map in \r{2.2} is a local diffeomorphism, we need to show $\partial_t \nu(t,x)\not=0$ in a conic neighborhood. Note that this map is well-defined as $\partial_x \phi(t,x)\not=0$. Assume that $\partial_t \nu(t,x)=0$. Then 
$$
	\frac{\partial_t \partial_x \phi(t,x)}{|\partial_x \phi(t,x)|} = \partial_x \phi(t,x)\frac{\partial_x \phi(t,x)\cdot \partial_t\partial_x \phi(t,x)}{|\partial_x \phi(t,x)|^3}
$$ 
	which implies that 
$$
	\partial_t \partial_x \phi(t,x)=c\partial_x \phi(t,x), \quad c= \frac{\partial_x \phi(t,x)\cdot \partial_t\partial_x \phi(t,x)}{|\partial_x \phi(t,x)|^2}.
$$
	This contradicts with the fact that $\partial_x \phi(t,x)$ and $\partial_t \partial_x \phi(t,x)$ are linearly independent. Now by Inverse Function Theorem, the map \r{2.2} is a local diffeomorphism as it is smooth and its Jacobian is nowhere vanishing.
\end{proof}
	One can extend the function $\phi$ to a homogeneous function of order one as follow:
\be{3.2}
	\varphi(x,\theta)=\psi_{\arg \theta} ^{-1}(x)\cdot \theta=|\theta|\phi(\arg \theta,x), \quad \text{where $\theta=(\theta^1, \theta^2)=|\theta|(\cos t, \sin t)\in \R^2 \setminus 0$}.
\ee 
	As we pointed out above, we work locally in a conic neighborhood of $t_0$ and $s_0$. This guarantees that function $\arg \theta$ is single-valued. To connect the local Bolker condition to Fourier Integral Operator theory, we have the following proposition.
\begin{proposition}
	For the function $\varphi$ defined by $\phi$ in \r{3.2}, the local Bolker condition \r{2.3} holds if and only if 
$$	
	\det \left( {\begin{array}{cc} \frac{\partial^2 \varphi}{\partial \theta^i \partial x^j} \end{array} } \right)\not = 0.
$$
\end{proposition}
\begin{proof}
	Since $\partial_{x} \varphi=|\theta|\partial_{x} \phi \not = 0$, we have
$$
	\frac{\partial^2 \varphi}{\partial \theta^1 \partial x^j}= \frac{\partial}{\partial \theta^1}(|\theta|\frac{\partial \phi}{ \partial x^j})= \frac{\theta^1}{|\theta|} \frac{\partial \phi}{ \partial x^j} - \frac{\theta^2}{|\theta|} \frac{\partial^2 \phi}{\partial t \partial x^j},
$$
	and 
$$
	\frac{\partial^2 \varphi}{\partial \theta^2 \partial x^j}= \frac{\partial}{\partial \theta^2}(|\theta|\frac{\partial \phi}{ \partial x^j})= \frac{\theta^2}{|\theta|} \frac{\partial \phi}{ \partial x^j} + \frac{\theta^1}{|\theta|} \frac{\partial^2 \phi}{\partial t \partial x^j},
$$
	where $t=\arg \theta$. Assume first that $\partial_x \phi(t,x)$ and $\partial_t \partial_x \phi(t,x)$ are linearly independent. We show that columns in the matrix $\frac{\partial^2 \varphi}{\partial \theta^i \partial x^j}$ are linearly independent for $i=1,2$. So let
$$
	c_1\frac{\partial^2 \varphi}{\partial \theta^1 \partial x^j} + c_2\frac{\partial^2 \varphi}{\partial \theta^2 \partial x^j}=0.
$$
	Then we have
$$
	(c_1\frac{\theta^1}{|\theta|} + c_2\frac{\theta^2}{|\theta|})\frac{\partial \phi}{ \partial x^j} +( -c_1 \frac{\theta^2}{|\theta|} +  c_2\frac{\theta^1}{|\theta|}) \frac{\partial^2 \phi}{\partial t \partial x^j}=0.
$$
	Since $\partial_x \phi(t,x)$ and $\partial_t \partial_x \phi(t,x)$ are linearly independent, we have
$$
	c_1\theta^1 + c_2\theta^2=0 , \quad -c_1 \theta^2 +  c_2\theta^1=0, 
$$
	which simply implies that $c_1=c_2=0$, and therefore $\frac{\partial^2 \varphi}{\partial \theta^i \partial x^j}$ are linearly independent for $i=1,2$.
	
	Assume now that $\frac{\partial^2 \varphi}{\partial \theta^i \partial x^j}$ are linearly independent for $i=1,2$. We show that $\partial_x \phi(t,x)$ and $\partial_t \partial_x \phi(t,x)$ are linearly independent. We first rewrite $\partial_x \phi(t,x)$ and $\partial_t \partial_x \phi(t,x)$ as follow:  
$$
	\theta^1\frac{\partial^2 \varphi}{\partial \theta^1 \partial x^j}=\frac{(\theta^1)^2}{|\theta|} \frac{\partial \phi}{ \partial x^j} - \frac{\theta^1 \theta^2}{|\theta|} \frac{\partial^2 \phi}{\partial t \partial x^j}, 
$$
	and
$$
	\theta^2\frac{\partial^2 \varphi}{\partial \theta^2 \partial x^j}=  \frac{(\theta^2)^2}{|\theta|} \frac{\partial \phi}{ \partial x^j} + \frac{\theta^1\theta^2}{|\theta|} \frac{\partial^2 \phi}{\partial t \partial x^j}.
$$
	Adding the last two equations we get
$$
	\frac{\theta^1}{|\theta|}\frac{\partial^2 \varphi}{\partial \theta^1 \partial x^j} + \frac{\theta^2}{|\theta|}\frac{\partial^2 \varphi}{\partial \theta^2 \partial x^j}= \frac{\partial \phi}{\partial x^j}.
$$
	Consider 
$$
	-\theta^2\frac{\partial^2 \varphi}{\partial \theta^1 \partial x^j}=-\frac{\theta^1\theta^2}{|\theta|} \frac{\partial \phi}{ \partial x^j} + \frac{(\theta^2)^2}{|\theta|} \frac{\partial^2 \phi}{\partial t \partial x^j}, 
$$
	and
$$
	\theta^1\frac{\partial^2 \varphi}{\partial \theta^2 \partial x^j}=  \frac{\theta^1\theta^2}{|\theta|} \frac{\partial \phi}{ \partial x^j} + \frac{(\theta^1)^2}{|\theta|} \frac{\partial^2 \phi}{\partial t \partial x^j}.
$$
	Adding the last two equations, we have
$$
	-\frac{\theta^2}{|\theta|}\frac{\partial^2 \varphi}{\partial \theta^1 \partial x^j} + \frac{\theta^1}{|\theta|}\frac{\partial^2 \varphi}{\partial \theta^2 \partial x^j}=  \frac{\partial^2 \phi}{\partial t \partial x^j}.
$$
	Now assume that  
$$
	\tilde{c}_1\frac{\partial \phi}{\partial x^j} + \tilde{c}_2\frac{\partial^2 \phi}{\partial t \partial x^j}=0.
$$ 
	In a similar way as we showed above and using the fact that $\frac{\partial^2 \varphi}{\partial \theta^i \partial x^j}$ are linearly independent for $i=1,2$, we conclude that $\tilde{c}_1=\tilde{c}_2=0$. This proves the proposition. 
\end{proof}
	In principle, Proposition 2.2 implies that we can use our analysis with \r{2.3}, see \cite{1}. 
\begin{definition}
	We say that $(x_0,\xi^0)\in \R^n \times (\R^n \setminus 0)$ is not in the Wave Front Set of $f \in \mathcal{D}^\prime (\R^n)$, $\WF f$, if there exists $\phi \in C^\infty _0(\R^n)$ with $\phi(x_0)\not = 0$ so that for any $N$, there	exists $C_N$ such that
$$
	|\hat{\phi f}(\xi)| \le C_N(1+|\xi|)^{-N}
$$
	for $\xi$ in some conic neighborhood of $\xi^0$. 
\begin{remark}
		The above definition is independent of the choice of $\phi$.
\end{remark}
\end{definition}
\begin{definition}
	For the case of a scalar-valued distribution, define the Analytic Wave Front Set, $\WFA (f)$, as the complement of all $(x,\xi) \in T^*(\R^n\setminus 0)$ such that
$$
	\int e^{\i\lambda |x-y|\cdot\xi -{\lambda\over2}|x-y|^2}\chi(y)f(y)\d y = \mathcal{O}(e^{-\lambda\over C}), \quad \lambda>0
$$
	with some $C > 0$ and $\chi \in C^\infty _0$ equal to 1 near $x$.
\begin{remark}
	We recall that, there are three equivalent definitions of Analytic Wave Front Set in the literature due to Bros-Iagolnitzer \cite{4}, H\"{o}rmander \cite{16}, and Sato \cite{27}. Bony \cite{3} and Sj\"{o}strand \cite{28} have shown the equivalence of all these definitions.
\end{remark} 
\end{definition}
%\begin{definition}
%	We call $(x,\xi)\in T^*X\setminus 0$ a position singularity if $(x,\xi)\in \WF f$. We call $(x,\xi)\in T^*X\setminus 0$ a visible singularity if $(x,\xi)\in \WF f$ satisfies the visibility condition by Definition 2.2. Similarly, we call $(s,t,\sigma,\tau)$ a measurement singularity if $(s,t,\sigma,\tau)\in \WF \mathcal{Af}$.
%\end{definition}	
\section{Microlocal Analyticity}
	In this section, we study the microlocal analyticity of operator $\mathcal{A}$ for a given $f$.  We first compute the adjoint operator.
	
	{\it{Adjoint Operator $\mathcal{A^*}$.}} Let $\phi \in C^{\infty}(\R \times \tilde{X})$ be given, where $X$ is embedded in an open set $\tilde{X}$. We extend our function $f$ to be zero on $\tilde{X}\setminus X$. Consider now the one-dimensional level curves
$$
	H(s,t)=\lbrace x\in \tilde{X} : s=\phi(t,x) \rbrace, \quad s\in \R,\ t\in \R
$$
	with Euclidean measure $\d S_x$ induced by the volume form $\d x$ in the domain $X$. There exists a non-vanishing and smooth function $J(t,x)$ such that 
$$
	\d S_{s,t}(x) \wedge \d s = J(t,x) \d x.
$$ 
	Therefore,
$$	
	\int_{T_1} ^{T_2}\int_{\R} (\mathcal{A}f)\bar{g} \d s \d t=\int_{T_1} ^{T_2}\int_{\R}\int_{H(s,t)} \mu(t,x)f(x)\bar{g}(s,t) \d S_{s,t} \d s \d t$$ $$
	= \int_{T_1} ^{T_2}\int_{\tilde{X}} \mu(t,x)f(x)\bar{g}(\phi(t,x),t)J(t,x) \d x \d t,
$$
	where $T_1<t<T_2$ and $0<T_2-T_1\ll 2\pi$. In the second equality above, we used the fact that the double integral $\int_{\R}\int_{H(s,t)}$ equals to an integral over $\tilde{X}$, by Fubini's Theorem. Thus, the adjoint of $\mathcal{A}$ in $L^2(X,\d x)$ is 
$$
	\mathcal{A}^*g(x)=\int_{\R} \bar{\mu}(t,x)\bar{J}(t,x)g(\phi(t,x),t)\d t,
$$ 
	where $\mu$ is supported in $\lbrace t\in \R: T_1<t<T_2 \rbrace$. In fact, the adjoint $\mathcal{A}^*g(x)$ is localized in $t$ and is an average over all lines or curves $H(s,t)$ that go through $x$.\\
	
	{\it{Schwartz Kernel.}} Now we compute the Schwartz kernel of the operator $\mathcal{A}$.
\begin{lemma}
	The Schwartz kernel $K_{\mathcal{A}}$ of $\mathcal{A}$ is 
$$
	K_{\mathcal{A}}(s,t,y)=\delta(s-\phi(t,y))\mu(t,y)J(t,y),
$$
	where $J(t,y)=|\d_y \phi|=(\sum|\partial_{y^j} \phi|^2)^{\frac{1}{2}}$.
\end{lemma}
\begin{proof}
	Let $\Phi(s,t,y)=s-\phi(t,y)$. By \r{2.1} we have
$$
	\mathcal{A} f(s,t)=\int_{\phi(t,y)=s} {\mu}(t,y) f(y) \d S_{s,t}
	= \int_{\phi(t,y)=s} {\mu}(t,y) f(y) |\d_y \Phi| |\d_y \Phi|^{-1}\d S_{s,t}.
$$
	Since $\partial_{y^j} \Phi=-\partial_{y^j} \phi$ and $\partial_{y^j} \phi\not = 0$ when $\Phi=0$, by Theorem (6.1.5) H\"{o}rmander \cite{14}, we have 
$$
	|\d_y \Phi|^{-1}\d S_{s,t}=\Phi^*\delta_0.
$$
	Here $^*$ is pullback with $\Phi^*\delta_0=\delta_0 \circ \Phi$. 
	The second integral above can be written as
$$
	\int \Phi^*\delta_0  {\mu}(t,y) f(y) |\d_y \Phi| \d y = \langle \Phi^*\delta_0 {\mu} |\d_y \Phi|,f\rangle.
$$
	Therefore, the Schwartz kernel of $\mathcal{A}$ is
$$	
	K_{\mathcal{A}}(s,t,y)=\delta(s-\phi(t,y))\mu(t,y)|\d_y \Phi|.
$$ 
\end{proof} 
\begin{remark}
	One can compute the Schwartz kernel of $\mathcal{A}^*$ and $\mathcal{N}= \mathcal{A^*}\mathcal{A}$:
$$	
	K_{\mathcal{A}^*}(s,t,x)=\delta(\phi(t,x)-s)\bar{\mu}(t,x)J(t,x),
$$
$$	
	K_{\mathcal{N}}(s,t,x,y)= \int_\R \delta(\phi(t,x)-\phi(t,y)) \bar{\mu}(t,x)J(t,x)\mu(t,y)J(t,y)\d t.
$$ 
\end{remark}	
	The following lemma shows that the operator $\mathcal{A}$ is an elliptic Fourier Integral Operator (FIO).
\begin{lemma}
	Let $M=\lbrace (s,t,x): \Phi(s,t,x)=s-\phi(t,x)=0 \rbrace \subset Y \times X.$ Then the operator $\mathcal{A}$ is an elliptic FIO of order $-\frac{1}{2}$ associated with the conormal bundle of $M$:
$$
	N^* M=\lbrace(s,t,x,\sigma,\tau,\xi)\in T^*(Y\times X) \big | \quad (\sigma,\tau,\xi)=0 \quad \textbf{\normalfont on} \quad T_{(s,t,x)} M\rbrace,
$$ 
	where $(s,t,\sigma,\tau)$ and $(x,\xi)$ are the coordinates on $T^*Y$ and $T^*X$, respectively.
\end{lemma}

\begin{proof}
	By Lemma 3.1 the Schwartz kernel $K_{\mathcal{A}}$ has singularities conormal to the manifold $M$. Since $\dim X = \dim Y=2$, the Schwartz kernel $K_{\mathcal{A}}$ is conormal type in the class $I^{-\frac{1}{2}} (Y\times X;M)$, see (Section~18.2, \cite{14}). This shows that the operator $\mathcal{A}$ is an elliptic FIO of order $-\frac{1}{2}$ associated with the conormal bundle $N^* M.$ Note that $\sigma$ is a one-dimensional non-zero variable. 
\end{proof}
	We now compute the canonical relation $\mathcal{C}$ and show it is a four-dimensional non-degenerated conic submanifold of $N^* M$ parametrized by $(t,x,\sigma)$. Note that $N^* M$ is a Lagrangian submanifold of $T^*(Y\times X)$.  
		
\begin{proposition}
	Let $\mathcal{C}$ be the canonical relation associated with $M$. Then  
$$
	\mathcal{C}=\lbrace (\phi(t,x),t, \sigma, -\sigma\partial_t\phi(t,x); x, \sigma \partial_x\phi(t,x) \big | (\phi(t,x),t,x)\in M, \quad 0\not = \sigma \in \R \rbrace.
$$ 
	Furthermore, the canonical relation $\mathcal{C}$ is a local canonical graph if and only if for any $t,$ the map
\be{4.1}
	x \rightarrow \left( \phi(t,x),  \partial_t \phi(t,x) \right)
\ee
	is locally injective and local Bolker condition \r{2.3} holds.\\
\end{proposition}
\begin{proof}
	The twisted conormal bundle of M:
$$
	\mathcal{C} = (N^* M\setminus{0})^\prime=\lbrace(s,t,\sigma,\tau;x,\xi) \big | \quad (s,t,\sigma,\tau;x,-\xi)\in N^* M\rbrace,
$$
	gives the canonical relation associated with $M$. We first calculate the differential of the function $\Phi(s,t,x)=s-\phi(t,x)$. We have
$$
	\d\Phi(s,t,x)=\d s-\partial_t\phi(t,x)\d t - \partial_x \phi(t,x) \d x.
$$
	Therefore, the canonical relation is given by 
$$ 
	\mathcal{C}=\lbrace (\phi(t,x),t, \sigma, -\sigma\partial_t\phi(t,x); x, \sigma \partial_x\phi(t,x) \big | (\phi(t,x),t,x)\in M, \quad 0\not = \sigma \in \R \rbrace.
$$
	Now consider the microlocal version of double fibration:
	\begin{center}
		\begin{tikzcd}[]
			&  \mathcal{C} \arrow{dl}[swap]{\Pi_{Y}}\arrow{dr}{\Pi_X} & \\
			T^*(Y) && T^*(X)
		\end{tikzcd}
	\end{center}
	where	
$$
	\Pi_X(\phi(t,x),t,\sigma,-\sigma \partial_t \phi;x,\sigma \partial_x \phi)=(x,\sigma \partial_x \phi),$$ $$\Pi_Y(\phi(t,x),t,\sigma,-\sigma \partial_t \phi;x,\sigma \partial_x \phi)=(\phi(t,x),t,\sigma,-\sigma \partial_t \phi).
$$
	Our goal is to find out when the Bolker condition (locally) holds for $\mathcal{C}$, that is, $\Pi_Y:\mathcal{C}\rightarrow T^*(Y)$ is an injective immersion. We first compute its differential:
$$
	\d_{t,x,\sigma}\Pi_Y=\begin{pmatrix}
	\partial_t \phi & \partial_{x^1} \phi & \partial_{x^2} \phi & 0 \\
	1 & 0 & 0 & 0 \\ 
	0 &0 & 0 & 1\\
	-\sigma \partial^2 _t \phi & -\sigma \partial^2 _{t, x^1} \phi & -\sigma \partial^2 _{t, x^2} \phi &  \partial_t \phi
	\end{pmatrix}.
$$
	If $\d_{t,x,\sigma}\Pi_Y$ has rank equal to four, then the Bolker condition is locally satisfied. Indeed, this is true, as $\d_{t,x,\sigma}\Pi_Y$ has rank equal to four if and only if the condition \r{2.3} holds. This implies that $\dim \mathcal{C}=4$. Since the map in \r{4.1} is one-to-one, the projection $\Pi_Y: \mathcal{C}\rightarrow T^*(Y)$ is an injective immersion. Hence, $\Pi_Y$ is a local diffeomorphism. 
\end{proof}	
	The following lemma states whether position singularities and measurement singularities can affect each other. We refer the reader to Definitions 3.3 and 3.4, for position and measurement singularities.
\begin{lemma}
	Let $X$ be a fixed open set in $\R^2$ and $Y$ be the open sets of lines determined by $(s,t)$ in $\R^2$. Then, the map 
$$	
	\Pi_X \circ \Pi^{-1} _Y: T^*(Y) \longrightarrow T^*(X) 
$$ 
	is a local diffeomorphism.
\end{lemma}
\begin{proof}
	Consider the map $\Pi_Y:\mathcal{C}\rightarrow T^*(Y)$. We show that for a given $(s,t,\sigma,-\sigma \partial_t \phi)\in T^*(Y)$, one can determine $(x,\xi)\in T^*(X)$. Since $\partial_{t}\phi$ is non-zero ($\sigma$ and $\sigma \partial_{t}\phi$ are both non-zero,) for a given $(s,t)$ there exists a tangent vector to each level curve $H(s,t)$. By Remark 2.1, one can find a non-zero normal vector $\partial_{x}\phi$ on each level curve, and therefore $\xi=\sigma \partial_{x}\phi$. On each level curve $H(s,t)$, we have $s=\phi(t,x)$. Since $\partial_{x}\phi\not=0$, the Implicit Function Theorem implies that the variable $t$ determines $x$. Hence, the map $\Pi_Y$ is a local diffeomorphism. 
	
	Now consider the map $\Pi_X:\mathcal{C}\rightarrow T^*(X)$. Our goal is to determine $(s,t,\sigma,-\sigma \partial_t \phi) \in T^*(Y)$, for a given $(x,\xi)=(x,\sigma \partial_x \phi) \in T^*(X)$. By Proposition 2.1, the map
$$
	\frac{\xi}{|\xi|}=\frac{\partial_x \phi}{|\partial_x \phi|} \longrightarrow t,
$$ 
	is a local diffeomorphism for a fixed point $x$ provided that the condition \r{2.3} holds. Thus, $(x,\frac{\xi}{|\xi|})$ determines the variable $t$. In particular, for a given $(x,\xi)$ this implies that one can identify the level curve $H(s,t)$, as $(t,x)$ determines $\phi$, and therefore $s$ (on each level curve we have $s=\phi(t,x)$.) Since $\xi=\sigma \partial_x \phi$ with $\xi\not=0$, one can determine $\sigma=\frac{|\xi|}{|\partial_{x}\phi|}$. To determine the last variable $\sigma \partial_{t}\phi$, it is enough to take the partial derivative of $\phi$ with respect to the variable $t$. Thus, the map $\Pi_X$ is a local diffeomorphism. We remind that the above argument is valid when the condition \r{2.3} is satisfied. 
		
	Now since dim$(Y)$=dim$(X)$ and $\Pi_X: \mathcal{C}\rightarrow T^*(X)$ and $\Pi_Y: \mathcal{C}\rightarrow T^*(Y)$ are local diffeomorphisms, the map 
$$
	\Pi_X \circ \Pi^{-1} _Y: T^*(Y) \longrightarrow T^*(X) 
$$
$$
	(s,t,\sigma,\tau) \longmapsto (x,\xi) 
$$
	will be a local diffeomorphism.
\end{proof}	

\begin{remark}
	i) Note that, by Proposition 4.1.4 (H\"{o}rmander \cite{15}), if we show one of the maps $\Pi_Y$ or $\Pi_X$ is a local diffeomorphism, then the other map is also a local diffeomorphism as dim$(Y)$=dim$(X)$. We, however, in above lemma have shown that both maps are local diffeomorphisms, as the proof reveals whether each map will be a global diffeomorphism or not. In fact, for a fixed $(x,\xi)\in T^*(X)$ there might be more than one curve which resolves the same singularity.\vspace{1.5mm}
	\\
	ii) The map $\Pi_Y:\mathcal{C}\rightarrow T^*(Y)$ being a local diffeomorphism implies that one can always track the position singularities $(x,\sigma \partial_x \phi)\in \WF(f)$ by having the measurement singularities $(\phi(t,x),t,\sigma,-\sigma \partial_t \phi$ $(x,\xi))\in \WF (\mathcal{A}f)$.\vspace{1.5mm}
	\\
	iii) From the geometrical point of view, the map $\Pi_X:\mathcal{C}\rightarrow T^*(X)$ being a local diffeomorphism means that for any fixed position $x$ and covector $\xi$, there exists a curve (NOT necessarily unique) passing through $x$ perpendicular to $\xi.$ This means singularities in data, i.e. $(x,\sigma \partial_x \phi)\in \WF(f)$, can affect the measurement singularities, i.e. $(\phi(t,x),t,\sigma,-\sigma \partial_t \phi$ $(x,\xi))\in \WF (\mathcal{A}f)$. \vspace{1.5mm}
	\\
	iv) Proposition 4.1 and Lemma 4.3 show the local surjectivity of the map 
$$
	[T_1,T_2]\ni t\rightarrow \frac{\partial_x \phi(t,x)}{|\partial_x \phi(t,x)|} \in S^1, \quad  \text{ for a fixed $x$.}
$$
	Note that if the visibility condition holds, then we have the global surjectivity on $S^1.$
\end{remark}
	
\section{Global Bolker Condition} 
	In this section, we study the microlocalized version of the normal operator $\mathcal{N}=\mathcal{A}^*\mathcal{A}$ to prove a stability estimate. It is known that the normal operator $\mathcal{N}$ is a $\Psi$DO if %the global Bolker condition
	the projection $\Pi_Y: \mathcal{C}\rightarrow T^*(Y)$ is an injective immersion (see Proposition ~8.2, \cite{8}). For our analysis, in addition to the visibility and local Bolker conditions, we assume that the semi-global Bolker condition is satisfied which is similar to the {\it{No Conjugate Points}} assumption for the geodesics ray transform studied in (\cite{7,24}). 
		
	We first perform the microlocalization in constructing the operator $\mathcal{N}$ in a small conic neighborhood of a fixed covector $(x_0,\xi^0) \in T^*X \setminus 0$. By the visibility condition, there exists some $(s_0,t_0)$ such that $\phi(t_0,x_0)=s_0$ and $\partial_{x}\phi(t_0,x_0)\parallel \xi^0$; which means for each point $x_0$ and co-direction $\xi^0$, there exists a curve passing through $x_0$ where $\xi^0$ is normal to it. By semi-global Bolker condition, there exists a pair of neighborhoods of $(x_0,(s_0,t_0))$, $V$ and $U$, such that for any $(x,(s,t))\in V\times U,$ the visibility condition is preserved under small perturbations in $t$ variable. We now shrink $V$ and $U$ sufficient enough such that the local Bolker condition is also satisfied.   
	
	Define $\mathcal{N}= \chi_{_X}\mathcal{A^*}\chi_{_Y}\mathcal{A},$ where $\chi_{_X} (x)$ and $\chi_{_Y} (s,t)$ are non-negative cut-off functions in a neighborhood of $x_0$ and $(s_0,t_0)$, respectively, with property that the projections $\Pi_X: \mathcal{C}\rightarrow T^*(X)$ and $\Pi_Y: \mathcal{C}\rightarrow T^*(Y)$ are embeddings above $\supp(\chi_{_X})$ and $\supp(\chi_{_Y})$. In fact, the smooth cut-off functions $\chi_{_X}$ and $\chi_{_Y}$ are localizations on the base variables $x$ and $(s,t)$ and they are not $\Psi$DOs. The following theorem shows that the (microlocalized) normal operator $\mathcal{N}= \chi_{_X}\mathcal{A^*}\chi_{_Y}\mathcal{A}$ is a $\Psi$DO of order $-1$. 	
	
%	\begin{remark}
%		There might be some points $(s,t)$ which do not satisfy the visibility, semi-global and local Bolker conditions. A point for which all above three conditions are satisfied is called a $\bf{Regular}$ point. 
%	\end{remark}
	
\begin{theorem}
	Let $(x_0,\xi^0) \in T^*X \setminus 0$ be a fixed covector. Assume that the visibility, the local and semi-global Bolker conditions are satisfied near $(x_0,\xi^0)$. Let $\chi_{_X} $ and $\chi_{_Y}$ be non-negative cut-off functions defined above. Then the operator $\mathcal{N}= \chi_{_X}\mathcal{A^*}\chi_{_Y}\mathcal{A}$ is a classical $\Psi$DO of order $-1$ with principal symbol 
$$
	p(x,\xi)=(2\pi)^{-1}\chi_{_X}\frac{W(x,x,\xi)+W(x,x,-\xi)}{\tilde{h}(x,\xi)}
$$
	near $(x_0,\xi^0)$. 
	%the local diffeomorphism $\Pi_X \circ \Pi^{-1} _Y$) . 
	The functions $W$ and $\tilde{h}$ are defined as  
$$
	W(x,x,\xi)=\chi_{_Y}(\phi(t,x),t)|\mu(t,x)|^2 J^2(t,x), \quad \textbf{\normalfont and} \ \quad \tilde{h}(x,\xi)=\frac{|\xi|}{|\partial_x \phi(t,x)|}h(t,x),
$$
 	where $t=t(x,\xi)$ is well-defined locally by Lemma 4.3. 
\end{theorem} 
\begin{proof}
	For the proof we mainly follow (Lemma~ 2, \cite{17}). By the equation \r{2.1}, we have
$$
	\chi_{_Y}(s,t)\mathcal{A} f(s,t)=\int_{\phi(t,x)=s} \chi_{_Y}(\phi(t,x),t){\mu}(t,x) f(x) \d S_{s,t}.
$$
	Considering the Schwartz kernel of the microlocalized normal operator $\mathcal{N}=	\chi_{_X}(x)\mathcal{A}^*\chi_{_Y}(s,t)\mathcal{A}$, we split the integration over $\R$ into $\{\sigma>0\}$ and $\{\sigma<0\}$. We have
$$
	K_\mathcal{N}= \int_{\R}  \int_0 ^{+\infty}  e^{\i (\phi(t,x)-\phi(t,y))\sigma} \chi_{_X}(x)W(t,x,y) \d \sigma \d t
$$ 
$$
	+\int_{\R}  \int_0 ^{+\infty}  e^{-\i (\phi(t,x)-\phi(t,y))\sigma} \chi_{_X}(x)W(t,x,y) \d \sigma \d 	t=K_{\mathcal{N}^+}+K_{\mathcal{N}^-},
$$  
	where $K_{\mathcal{N}^+}$ and $K_{\mathcal{N}^-}$ are the Schwartz kernels of the operators $\mathcal{N}^+$ and $\mathcal{N}^-$ with $\mathcal{N}=\mathcal{N}^+ + \mathcal{N}^-$. We first consider $K_{\mathcal{N}^+}$. Note that $K_{\mathcal{N}^+}$, localized as the function $\phi$, priori satisfies the local Bolker condition \r{2.3}. By semi-global Bolker condition \r{2.4}, we have
$$
	\left\{\begin{array}{ll}
	\phi (t, x)=\phi (t, y)=s\\
	\partial_t\phi (t, x)=\partial_t\phi (t, y)\\
	\end{array}
	\right. \ \Longrightarrow \ x=y.
$$
	Now a stationary phase method implies that $K_{\mathcal{N}^+}$ is smooth away from the diagonal $\{x=y\}$. Since $\partial_x\phi (t, x)\not =0$, for a fixed $x$ there exists a neighborhood $\mathcal{U}$ on which we have normal vectors. We work on normal coordinates $(x^i, y^i)$ as coordinates on $\mathcal{U}\times \mathcal{U}$, with $x^i=y^i$. In these local coordinates, one can expand the phase function near the diagonal $\{x=y\}$. Let
\be{5.1}
	(\phi(t,x)-\phi(t,y))\sigma = (x-y) \cdot \xi(t,\sigma,x,y),
\ee
	where $\xi(t,\sigma,x,y)$ is defined by the map
$$		 
	(t,\sigma)\rightarrow\xi(t,\sigma,x,y)=\int_0 ^1 \sigma\partial_x\phi(t, x+\tau(y-x))\d \tau.
$$ 
	On the diagonal, we have $\xi(t,\sigma,x,x)= \sigma\partial_x\phi(t, x)=\xi$ and the map is a smooth diffeomorphism as
$$
	\det (\frac{\partial \xi}{\partial (t,\sigma)})_{\big |_{x=y}} = \det \left( {\begin{array}{cc}
	\frac{\partial \phi}{ \partial x^j},  \sigma\frac{\partial^2 \phi}{\partial t \partial x^j}       \end{array} } \right)=\sigma h(t,x)\not=0.
$$
	Notice that $\sigma=\frac{|\xi|}{|\partial_x \phi(t,x)|}$ and $t=t(x,\xi)$ is locally well-defined by Lemma 4.3. Therefore,
$$
	\tilde{h}(x,\xi)=\frac{|\xi|}{|\partial_x \phi|}h(t,x)\not=0.
$$
	Using the above change of variable $\r{5.1}$ on the diagonal yields 
$$
	K_{\mathcal{N}^+}(s,t,x,y)= (2\pi)^{-1} \iint_{\R^2}  e^{\i(x-y) \cdot \xi} \ \chi_{_X}(x) W(x,y,\xi) |\tilde{h}(x,\xi)|^{-1} \d \xi,
$$ 
	where the function $W$ is defined above. By restricting the amplitude to diagonal $\{x=y\},$ one can find the principal symbol of $K_{\mathcal{N}^+}$. Now the principal symbol of $K_{\mathcal{N}}$ is given by the sum of those for $K_{\mathcal{N}^+}$ and $K_{\mathcal{N}^-}$. Since the weight $\mu$ is a positive real analytic function, the normal operator $\mathcal{N}$ is a classical $\Psi$DO with principal symbol $p(x,\xi)$ provided the function $\phi$ satisfies the local and semi-global Bolker condition. Now since $\mu$ is nowhere vanishing and by local Bolker condition \r{2.3} $h(t,x)\not=0$, the operator $\mathcal{N}$ is an elliptic $\Psi$DO if the visibility condition is satisfied.  
\end{proof}
	
\section{Analysis of Global Problem and Stability}
	In previous sections, we studied the operators $\mathcal{A}$ and $\mathcal{N}$. We showed that under the visibility, local and semi-global Bolker conditions, the microlocalized normal operator $\mathcal{N}$ is a $\Psi$DO of order $-1$ in a small conic neighborhood of a fixed covector $(x,\xi) \in T^*X \setminus 0$. 
	
	To reconstruct $f\in L^2(X)$ from its measurements $\mathcal{A}f$ using the operator $\mathcal{N}$, we need to expand our results globally. As we pointed out in the begining of section five, the visibility, local and semi-global Bolker conditions (which are open conditions in a small conic neighborhood of $(x_0,\xi^0)$) are required for the analysis. We also employ non-negative cut-off functions $\chi_{_X}$ and $\chi_{_Y}$ in neighborhoods of $x_0$ and $(s_0,t_0)$, where the projections $\Pi_X$ and $\Pi_Y$ are embeddings above $\supp(\chi_{_X})$ and $\supp(\chi_{_Y})$. 
		
	Let $K\subset X$ be a compact subset and $(x_0,\xi^0)\in T^*K\setminus 0$ be a fixed covector. There exists a pair of conic neighborhoods $(\mathcal{V},\tilde{\mathcal{V}})$ with property $(x_0,\xi^0)\in \mathcal{V}$ and $\mathcal{V}\Subset \tilde{\mathcal{V}}$ such that the visibility, local and semi-global Bolker conditions are satisfied for $\tilde{\mathcal{V}}$. Let $\{\mathcal{V}_\alpha\}$ be an open covering for $T^*K\setminus 0$. Since $T^*K\setminus 0$ is conically compact subset of $T^*X\setminus 0$, by a compactness argument, there exists a finite subcover of $\{\mathcal{V}_i\}$. By Theorem 5.1, the microlocally restricted normal operators $\mathcal{N}_i=\chi_{i_X}\mathcal{A}^* \chi_{i_Y}\mathcal{A}$ are $\Psi$DOs of order $-1$ supported in a conic neighborhood $\mathcal{V}_i$ (where the visibility, local and semi-global Bolker conditions are satisfied), with the principal symbols
	$$
	p_{i}(x,\xi)=(2\pi)^{-1}\chi_{i_X}(x)\frac{W_i(x,x,\xi)+W_i(x,x,-\xi)}{\tilde{h}(x,\xi)},
	$$
	where  
	$$
	W_i(x,x,\xi)=\chi_{i_Y}(\phi(t,x),t)|\mu(t,x)|^2 J^2(t,x), \quad \ \quad \tilde{h}(x,\xi)=\frac{|\xi|}{|\partial_x \phi(t,x)|}h(t,x),
	$$
	and $t=t(x,\xi)$ is well-defined locally by Lemma 4.3. Here $\{\chi_{i_X}\}$ and $\{\chi_{i_Y}\}$ are families of smooth cut-off functions which are non-negative in  neighborhoods of $V_i\ni x_0$ and $U_i \ni (s_0,t_0)$, with property that $\supp \chi_{i_X} \subset {V}_i$ and $\supp \chi_{i_Y} \subset {U}_i$. We remind that, the smooth cut-off functions $\chi_{i_X}$ and $\chi_{i_Y}$ are localizations on the base variables $x$ and $(s,t)$ and they are not $\Psi$DOs.
	
	Set $\mathcal{N}=\sum \mathcal{N}_i.$ Now for any $(x,\xi)$, there exists $k$ such that $\chi_{k_X}(x)\not=0$ and all other terms are non-negative. Hence $\sum \mathcal{N}_i$ is elliptic, and therefore the operator $\mathcal{N}$ is a classical $\Psi$DO of order $-1$ with principal symbol $P(x,\xi)=\sum p_{i}(x,\xi)$.

\begin{remark}
	It should be noted that in our analysis, the cut-off functions are used for the $C^\infty$ results. For the case of analytic arguments, one cannot use cut-off functions.
\end{remark}
	
	In the following proposition, we show that for any neighborhood of a fixed covector $(x_0,\xi^0)\in T^* X \setminus 0$, ellipticity holds along normals in a conic neighborhood of this covector. We point out that, one can use the "eating away at $\supp f$" argument, first stated by Boman and Quinto \cite{2}, to conclude the similar results. 
\begin{proposition}
	Assume that the dynamic operator $\mathcal{A}$ satisfies the visibility, the local and semi-global Bolker conditions for all $(s,t)\in Y$ and $(x_0,\xi^0)\in T^*(X)\setminus 0 $. Let $\phi$ be a real analytic function and $\mu$ be a positive real analytic weight. Let $f\in L^2(X)$ with $\supp f \subset X$. If $\mathcal{A}f=0$ in a neighborhood of some level curves, $l_0$, determined by $(s_0,t_0)$, then 
$$
	\WFA(f)\cap N^*(l_0)=\emptyset.
$$
\end{proposition}
\begin{proof}
	Let $(x_0,\xi^0)\in T^* X \setminus 0$ be fixed. By the visibility condition, there exists $(s_0,t_0)$ such that $\phi(t_0,x_0)=s_0$ and $\partial_{x}\phi(t_0,x_0)\parallel \xi^0$. Now the proof follows directly from [Proposition~1, \cite{17}] and applying it to all conormals of the fixed curve $l_0,$ determined by $(s_0,t_0)$.
\end{proof}	
\begin{remark}
	For the results in Proposition 6.1, we only need the visibility, the local and semi-global Bolker conditions to be satisfied near $N^*(l_0).$ However, to conclude the following corollary, we need to have the above three conditions satisfied globally, i.e. for all $(s,t)\in Y$ and $(x_0,\xi^0)\in T^*(X)\setminus 0$.  
\end{remark}
\begin{corollary}
	Under the assumption of Proposition 6.1, $\mathcal{A}f=0$ implies that $f=0$.
\end{corollary}
\begin{proof}
	Let $\tilde{X} \supset \supp f$ be an open set where the function $f$ is extended to be zero on $\tilde{X}\setminus X$ ($X$ is embedded in the set $\tilde{X}$.) Consider all level curves intersecting $\tilde{X}.$ By visibility condition, there exists a level curve $l_0$ determined by $(s_0,t_0)$ such that $\phi(t_0,x_0)=s_0$ and $\partial_{x}\phi(t_0,x_0)\parallel \xi^0$ (i.e. each singularity is visible). On the other hand, the local and semi-global Bolker conditions guarantee that there exist some lines in the exterior of $\supp f$.
	% a neighborhood of $l_0$ such that for any line in that neighborhood 
	By assumption, $\mathcal{A}f=0$ for all these level curves. Now, Proposition 6.1 implies that $f$ is analytic in the interior of $\tilde{X}$. Since $f$ is identically zero on $\tilde{X}\setminus X$, $f$ must be identically zero on all of $X$. Hence $\mathcal{A}$ is injective.
\end{proof}	
	The following proposition is a standard stability estimate which follows from elliptic regularity see (Theorem~2, \cite{29}) and (Proposition ~V.3.1, \cite{30}). 
\begin{proposition}
	Let the real analytic function $\phi$ satisfies the visibility, the local and semi-global Bolker conditions and $\mu$ be a positive real analytic weight. Let $K$ be a compact subset of $X$. Then for all $f\in L^2(K)$ and $s>0$ there exists $C>0$ and $C_s>0$ depending on $s$ such that
$$
	\norm{f}_{L^2(K)} \leq C \norm {\mathcal{N}f}_{H^1(\tilde{X})} + 
	C_s \norm {f}_{H^{-s}}, \quad \text{$\forall s$.}
$$ 
	Moreover, if $\mathcal{N}:L^2(K)\rightarrow H^1(\tilde{X})$ is injective, then there exists a stability estimate,
$$
	\norm{f}_{L^2(K)} \leq C^\prime \norm {\mathcal{N}f}_{H^1(\tilde{X})}
$$
	where $C^\prime>0$ is a constant. 
	
\end{proposition}
\begin{proof}
	The proof directly follows from Theorem 5.1 and above arguments. 
\end{proof}
\begin{remark}
	Note that the way the parametrix is constructed in above proposition, one has control on how the constant $C$ to be chosen. This, however, is not the case for $C^\prime$ in the second inequality.  
\end{remark}
	In what follows, we perturb $\phi$ and $\mu$, and prove that the perturbation yields a small constant times an $L^2$-norm of the function $f$ which can be absorbed by the left-hand side of above estimate. The following lemma is in the spirit of [Lemma~4, \cite{17}].
\begin{lemma}
	Let $\mathcal{A}$ be a dynamic operator satisfying the visibility, the local and semi-global Bolker conditions with a real analytic function $\phi$ and positive real analytic weight $\mu$. There exists a $k \gg 2$ and $(\tilde{\phi},\tilde{\mu})\in C^k$ such that if 
$$
	\parallel\phi-\tilde{\phi}\parallel_{C^k(\R \times\tilde{X})}, \quad \norm {\mu-\tilde{\mu}}_{C^k(\R \times\tilde{X})} <\delta \ll 1,
$$
	then there exists $C\ge0$ depending on the $C^k(\R \times\tilde{X})$ norm of $\phi$ and $\mu$ such that 
$$
	\parallel(\mathcal{N}-\tilde{\mathcal{N}})f\parallel_{H^1(\tilde{X})} \leq C\delta \norm {f}_{L^2(\tilde{X})}.
$$
	Here 
$$	\mathcal{N}=\sum_i\mathcal{N}_i=\sum_i\chi_{i_X}\mathcal{A}^* \chi_{i_Y}\mathcal{A},
	\quad \quad \tilde{\mathcal{N}}=\sum_i\tilde{\mathcal{N}}_i=\sum_i\chi_{i_X}\tilde{\mathcal{A}^*} \chi_{i_Y}\tilde{\mathcal{A}}
$$ 
	are two microlocally restricted normal operators corresponding $\mu$ and $\tilde{\mu}$, respectively, and the cut-off functions $\chi_{i_X}$ and $\chi_{i_Y}$ are defined as above.
\begin{proof}
	Let $(x_0,\xi^0)\in T^*(X)\setminus 0$ be a fixed covector. By the visibility condition, there exists a line $l_0$, determined by $(s_0,t_0)$, such that $\phi (t_0,x_0) = s_0$ and $\partial_{x}\phi(t_0,x_0)\parallel\xi^0$. Let $\chi_{_X}$ and $\chi_{_Y}$ be smooth cut-off functions defined above in neighborhoods of $x_0$ and $l_0$ corresponding to $\phi\in C^k$ with $k$ large enough. By Lemma 4.3 and Remark 4.2, for any level curve $l$ close to $l_0$, a perturbation of $\phi\in C^k$ results in the perturbation of the family of the level curves near $\phi(t_0,x)=s_0$. Since the local and semi-global Bolker conditions are open conditions, the visibility condition is preserved under the small perturbation in a neighborhood of $l_0$. On the other hand, a priori, we assumed that $\phi$ and $\tilde{\phi}$ are $\delta$-close with $C^k$-topology. Therefore, one can choose the same cut-off function $\chi_{_X}$ and $\chi_{_Y}$ such that both projections $\Pi_Y$ and $\tilde{\Pi}_Y$ are embeddings on their support and the visibility, the local and semi-global Bolker conditions are satisfied in each neighborhood. Therefore for each $i$, Theorem 5.1 implies that the microlocally restricted normal operators $\mathcal{N}_i=\chi_{i_X}\mathcal{A}^* \chi_{i_Y}\mathcal{A}$ and $\tilde{\mathcal{N}}_i=\chi_{i_X}\tilde{\mathcal{A^*}}\chi_{i_Y}\tilde{\mathcal{A}}$ are elliptic $\Psi$DOs with symbols depending on $\phi$, $\mu$ and $\tilde{\phi}$, $\tilde{\mu}$, respectively. 
	
	We now directly apply the argument on [Lemma~4, \cite{17}] to $\mathcal{N}_i ^\pm-\tilde{\mathcal{N}}_i ^\pm$, to conclude that for each $i$
$$
	\parallel\mathcal{N}_i^\pm-\tilde{\mathcal{N}}_i^\pm\parallel_{L^2_c(\tilde{X})\rightarrow H^1(\tilde{X})} =\mathcal{O}(\delta),
$$	
	and hence,
$$
	\parallel(\mathcal{N}_i-\tilde{\mathcal{N}}_i)f\parallel_{H^1(\tilde{X})} \leq C\delta \parallel f \parallel_{L^2(\tilde{X})}.
$$
	Now the fact that the operator $\mathcal{N}$ is a finite sum of operators of the form $\mathcal{N}_i$, as well as using the triangle inequality 
$$
	\parallel(\mathcal{N}-\tilde{\mathcal{N}})f\parallel_{H^1(\tilde{X})} \leq \sum_i \parallel( \mathcal{N}_i -\tilde{\mathcal{N}}_i)f\parallel_{H^1(\tilde{X})},
$$
	conclude the results.
\end{proof}	
\end{lemma}
	Next result is a stability estimate for a generic class of dynamic operators satisfying the visibility, the local and semi-global Bolker conditions.
\begin{theorem}
	Let $X$ be an open set of points $($positions$)$ $x$ lying on lines in $Y$, where $Y$ is the open sets of lines determined by $(s,t)$ in $\R^2$. Let $\mathcal{A}:L^2(X)\rightarrow H^1(\tilde{X})$, satisfying the visibility, the local and semi-global Bolker conditions, be an injective dynamic operator defined by the real analytic function $\phi$ and positive real analytic weight $\mu$. Then\vspace{1.5mm} 
	\\
	i) For any $\tilde{\phi}\in \n(\phi)$ and $\tilde{\mu}\in \n(\mu)$ with $C^k$-topology $($$k$ an arbitrary large natural number$)$ and for all $f\in L^2(K)$ with $K$ a compact subset of $X$, there exists $C\ge 0$ such that 
$$
	\parallel f\parallel_{L^2(K)} \leq C \parallel \tilde{\mathcal{N}f} \parallel_{H^1(\tilde{X})}.
$$
	In particular, the operator $\tilde{\mathcal{A}}$ is injective.\vspace{1.5mm}
	\\
	ii) The following stability estimate remains true for any perturbation of $\phi$ and $\mu$: 
$$
	\parallel f\parallel_{L^2(K)}/C \ \leq \ \parallel \mathcal{N}f \parallel_{H^1(\tilde{X})} \ \leq \ C \parallel f\parallel_{L^2(K)}.
$$ 
\end{theorem}
\begin{proof}
	i) $\mathcal{A}$ is injective, thus by Proposition 6.2, we have the following stability estimate:
$$
	\parallel f \parallel_{L^2(\tilde{X})} \ \leq \ C_1 \parallel \mathcal{N}f \parallel_{H^1(\tilde{X})} \ = \ C_1 \parallel  \tilde{\mathcal{N}}f + (\mathcal{N} - \tilde{\mathcal{N}})f \parallel_{H^1(\tilde{X})}
$$
$$
	\leq \ C_1 \parallel \tilde{\mathcal{N}}f \parallel_{H^1(\tilde{X})} \ + \ C_1\parallel(\mathcal{N} - \tilde{\mathcal{N}})f \parallel_{H^1(\tilde{X})}.
$$
	By Lemma 6.1, there exists a constant $C_2\ge0$ such that
$$
	\parallel(\mathcal{N} - \tilde{\mathcal{N}})f \parallel_{H^1(\tilde{X})} \  \leq \  C_2\delta \parallel f\parallel_{L^2(\tilde{X})},
$$
	and therefore,
$$
	\parallel f \parallel_{L^2(\tilde{X})} \ \leq \ C_1\parallel \tilde{\mathcal{N}}f \parallel_{H^1(\tilde{X})} \  + \ C_1 C_2\delta \parallel f\parallel_{L^2(\tilde{X})}.
$$
	Letting $\delta < \min \{(2C_1C_2)^{-1},1/2\}$ yields
$$
	\parallel f \parallel_{L^2(K)} \  \leq \  C \parallel \tilde{\mathcal{N}}f \parallel_{H^1(\tilde{X})}.
$$ 
	Assume now that $\tilde{\mathcal{A}}f=0$. Then 
$$
	\tilde{\mathcal{N}}f =\sum_i\mathcal{\tilde{A}}^*\chi_i\tilde{\mathcal{A}}f = 0, \quad \text{as \ $\tilde{\mathcal{A}}f=0$.}
$$
	The last inequality above implies that $f=0$. Hence, the operator $\tilde{\mathcal{A}}$ is injective. \vspace{1.5mm}
	\\
	ii) This part follows directly from the first part and the continuity of pseudodifferential operator $\tilde{\mathcal{N}}$. 
\end{proof}
	
\begin{proof}[\textbf{Proof of Theorem 2.1}] 
		The proof directly follows from Theorem 6.1. 
\end{proof}

\section{Analysis of the Initial Dynamic Problem}
	In this section, we state the implications of our analysis for the partial case, where the dynamic operator is given by \r{1.1}. This
	%$$
	%\mathcal{A} f(s,t)=\int_{z\cdot \omega(t)=s} \mu(t,z) f(\psi_t (z)) \d S,
	%$$
	corresponds to the initial example of scanning the moving object with changing its shape. Some part of above results, the local and semi-global Bolker assumptions, are also given in [Theorem~14, \cite{11}] and the problem of recovery of singularities has been analyzed. The periodic and non-periodic motions with $\phi(t,x)=\psi_t ^{-1}(x)\cdot\omega(t)$ have been studied in \cite{11} to explain which singularities are visible (see Theorems 24, 26). 
		 
	Using a change of variable $x=\psi_t (z)$, the dynamic operator $\mathcal{A}$ can be written as:
	$$
	\mathcal{A} f(s,t)=\iint_{\R^2} J(t,x) \mu(t,\psi_t ^{-1}(x)) f(x)\ \delta(s-\psi_t ^{-1}(x)\cdot \omega(t)) \ dx, 
	$$ 
	where $\psi_t ^{-1}(x)\cdot \omega(t)$ is the level curve corresponding to $\mathcal{A}$.\vspace{1.5mm}
	\\
	{\it{Canonical relation.}} Setting $\Phi(s,t,x)=s-\psi_t ^{-1}(x)\cdot \omega(t)$ in Proposition 4.1, the canonical relation $\mathcal{C}$ associated with $\mathcal{A}$ will be
	$$
	\mathcal{C}=\lbrace (\psi_t ^{-1}(x)\cdot \omega(t),t, \sigma, -\sigma(\partial_t\psi_t ^{-1}(x)\cdot\omega(t)+\psi_t ^{-1}(x)\cdot \omega^{\bot}(t)); x, \sigma \partial_x \psi_t ^{-1}(x)\cdot\omega(t)) \big | (s,t,x)\in M\rbrace.
	$$
	The microlocal version of double fibration is given by:
	\begin{center}
		\begin{tikzcd}[]
			&  \mathcal{C} \arrow{dl}[swap]{\Pi_{Y}}\arrow{dr}{\Pi_X} & \\
			T^*(Y) && T^*(X)
		\end{tikzcd}
	\end{center}
	where
$$
	\Pi_X(\psi_t ^{-1}(x)\cdot \omega(t),t,\sigma,-\sigma\partial_t(\psi_t ^{-1}(x)\cdot\omega(t));x,\sigma \partial_x \psi_t ^{-1}(x)\cdot\omega(t))=(x,\sigma \partial_x \psi_t ^{-1}(x)\cdot\omega(t)),
$$
$$
	\Pi_Y(\psi_t ^{-1}(x)\cdot \omega(t),t,\sigma,-\sigma\partial_t(\psi_t ^{-1}(x)\cdot\omega(t));x,\sigma \partial_x \psi_t ^{-1}(x)\cdot\omega(t))=(\psi_t ^{-1}(x)\cdot \omega(t),t,\sigma,-\sigma\partial_t(\psi_t ^{-1}(x)\cdot\omega(t))).\vspace{1.5mm}
	\\
$$
	{\it{Visibility.}} The operator $\mathcal{A}$ satisfies in the visibility condition if for any $(x,\xi)\in T^*X\setminus 0$, the map given by \r{2.5} is locally surjective.\vspace{1.5mm}
	\\	
	{\it{Local Bolker Condition.}} 
	As it is shown in Proposition 4.1, the projection $\Pi_Y$ is an immersion if the matrix $\d_{t,x,\sigma}\Pi_Y$ has rank equal to four or equivalently $\det (\d_{t,x,\sigma}\Pi_Y)\not =0 $. Since
$$
	\det (\d_{t,x,\sigma}\Pi_Y )= \det \left( {\begin{array}{cc} 
	\frac{\partial \psi_t ^{-1}(x)\cdot\omega(t)}{ \partial x^j},  \frac{\partial^2 \psi_t ^{-1}(x)\cdot\omega(t)}{\partial t \partial x^j} \end{array} } \right) = h(t,x),
$$
	the projection $\Pi_Y$ being an immersion is equivalent to the condition \r{2.3} being non-zero, i.e. $h(t,x)\not=0.$ \vspace{1.5mm}
	\\	
	{\it{Semi-global Bolker condition (No conjugate points condition).}} By condition \r{2.4}, $\Pi_Y$ is injective if the map
$$
	x \rightarrow \left( \psi_t ^{-1}(x)\cdot\omega(t),  \partial_t (\psi_t ^{-1}(x)\cdot\omega(t)) \right)
$$
	is one-to-one. \vspace{1.5mm}
	\\		
	{\it{The normal operator $\mathcal{N}$ is a $\Psi$DO of order $-1$.}} Under the local and the semi-global Bolker conditions, Theorem 5.1 implies that the normal operator $\mathcal{N}$ associated with the dynamic operator $\mathcal{A}$ is a $\Psi$DO of order $-1$ with principal symbol $p(x,\xi)$ near each $(x_0,\xi^0)$. The principal symbol is given by 
$$
	p(x,\xi)=(2\pi)^{-1}|\partial_x \psi_t ^{-1}(x)\cdot\omega(t)|\frac{|\mu(x,\xi)|^2 J^2(x,\xi)+|\mu(x,-\xi)|^2 J^2(x,-\xi)}{|\xi|h(x,\xi)},
$$
	where $t=t(x,\xi)$ is locally well-defined by Lemma 4.3.
\begin{remark}
	Note that we do not require the function $\phi(t,x)$ to be smoothly periodic. 
\end{remark}

	\textbf{Fan Beam Geometry.} In previous sections, we showed that the dynamic operator $\mathcal{A}$ with $\phi(t,x)=\psi_t ^{-1}(x)\cdot\omega(t)$ in parallel beam geometry, belongs to a more general integral geometry problem. We formulated the visibility, local and semi-global Bolker conditions, and derived our results for the case when $\phi(t,x)=\psi_t ^{-1}(x)\cdot\omega(t)$.  
% Note also that the operator $\mathcal{A}$ can also be written in the following representation: 
%$$
%	\mathcal{A}_P f(s,\beta)=\int_{\R} f(s\omega(\beta) + p\omega^\bot (\beta)) \ dp. 
%$$
%	where $\omega(\beta)=(\cos \beta, \sin \beta).$
	
	Another common geometry which is often used in numerical simulations is Fan beam geometry. In this geometry, the assumption is that each scan is taken from a boundary point $S$ (Source) and all directions instantly but the object moves when we change $S$ (see Figure 1).
\begin{figure}[h!]
		\includegraphics[scale=0.6]{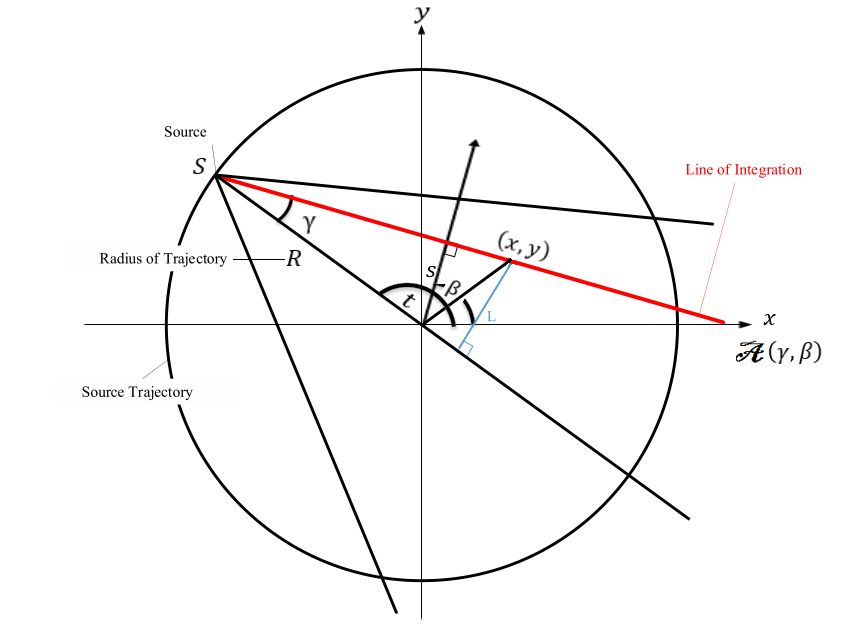}
		\caption{Parallel-Fan beam geometry relation.}
\end{figure} 

	Using the Parallel-Fan beam relation 
$$
	s=R\sin \gamma \ \ \  \quad    \quad \beta = t + \gamma - \frac{\pi}{2},
$$ 
	and finding an appropriate level curve $\phi$, one can show the dynamic operator $\mathcal{A}$ (in fan beam geometry) is also a special case of the general integral geometry problem discussed in this paper. Although the dynamic operator $\mathcal{A}$, will have different representations due to different time-parameterizations in parallel and fan beam geometries, they both can be categorized by the same general integral geometry problem.
	
	%	In this section, for a small enough motion,  our intention is to show that one can employ the perturbation arguments to conclude our general results for the dynamic operator $\mathcal{A}$ in fan beam geometry as the local and semi global Bolker conditions are open conditions, and they are clearly satisfied if there is no motion.  
	
	For simplicity, our analysis in this section is restricted to a static case, i.e. $\phi(t,x)=x\cdot\omega(t),$ as the visibility, the local and semi global Bolker conditions are clearly satisfied when there is no motion. One can achieve the same results for the general case where the motion is not necessarily small. We, however, do not provide details on how to formulate the visibility, the local and semi-global Bolker conditions and rather state that our results are valid if these conditions are satisfied.
	
	Let lines %(in fan beam geometry) 
	along which the dynamic operator of $f$ is known, are specified by $\gamma$ (the angle between the incident ray direction and the line from the source to the rotation center) and $t$ (the angular position of the source). Then the fan beam data at time $t$ is given by
$$
	\mathcal{A}_F f(t,\gamma)=\int^\infty _0 f(S(t) + p\theta(\gamma)) \ dp, \quad  \text{$\theta(\gamma)\in S^1$},
$$
	where $S(t)$ is the source at time $t$ which moves along the trajectory with radius $R$. Here $t$ is both a parameter along the source trajectory and the time variable. Note also that using the Parallel-Fan beam geometry relation one can derive the fan beam dynamic operator $\mathcal{A}_F$, given by 
$$
	\mathcal{A}_F f(t,\gamma)=\mathcal{A}_P f(R\sin \gamma, t + \gamma - \frac{\pi}{2}).
$$
	Since the Jacobian
$$ 
	\big |\frac{\partial(s,\beta)}{\partial(t,\gamma)}\big|= R\cos \gamma
$$
	is non-zero, the transformation between these two geometries is smooth.
	
	To implement our results in fan beam geometry, we need to find appropriate level curves $\phi$.	Let $S(t)$ be the source and $x$ be the point on the incident ray, see figure 1. We first set	
$$
	\boldsymbol{\vec{\alpha}}=x-S(t)=(x^1-R\cos t, x^2-R\sin t),  
$$
	and then compute the perpendicular vector $\boldsymbol{\vec{\alpha}}^\bot$ as follow:
$$
	\boldsymbol{\vec{\alpha}}^\bot=\frac{\sgn(x^1-R\cos t)}{|\boldsymbol{\vec{\alpha}}|} (R\sin t-x^2 , x^1-R\cos t)=(\cos \alpha^\bot, \sin \alpha^\bot).
$$
	We only work with one direction from two possible orientations for $\boldsymbol{\vec{\alpha}}^\bot$; say the one with $\boldsymbol{\sgn}(x^1-R\cos t)>0$. For a fixed point $x$ on the incident ray and a specific time $t$, the polar angle $\alpha^\bot$ is determined by
$$
	\alpha^\bot=\arg(\boldsymbol{\vec{\alpha}}^\bot)=\tan^{-1}\big(\frac{x^1-R\cos t}{R\sin t-x^2}\big).
$$
	We set 
$$
	\phi(t,x)=\arg(\boldsymbol{\vec{\alpha}}^\bot).
$$	
	Now our results are valid if the visibility, the local and semi-global Bolker conditions are satisfied for this choice of function $\phi$. Note that, here the function $\boldsymbol{\arg}$ is not globally defined but this does not affect the analysis, as our results are local and we have chosen the branch where $\boldsymbol{\sgn}(x^1-R\cos t)>0$. One can choose another branch of $\tan^{-1}$, however, this plays no role in differentiation which is involved in all above main three conditions.
	
	%\begin{remark}
	%	What we discussed above is the fan beam geometry for the case where there is no motion or the motion is small enough, i.e. $\psi_t\cong\Id$. For the general case, one can replace $x=(x^1,x^2)$ by  $\psi_t ^{-1}(x)=(\psi_{t} ^{-1}(x^1),\psi_{t} ^{-1}(x^2))$ to get the main result. 
	%\end{remark} 
	
\textbf{Acknowledgments.} The author would like to express his special gratitudes for Professor Plamen D. Stefanov for introducing the problem and his valuable discussions throughout this work. The author thanks Professor Todd Quinto for his helpful comments. The author also thanks referees for their valuable comments have helped in improving the manuscript.


\begin{thebibliography}{1}
	
	\bibitem{1}
	G. Beylkin.
	\textit{The inversion problem and applications of the generalized Radon transform.}
	Comm. Pure Appl. Math., 37 (1984), pp. 579–599.
	
	\bibitem{2}
	J. Boman and E. T. Quinto.
	\textit{Support theorems for real-analytic Radon transforms.} 
	Duke Math. J., 55 (1987), 943–948.
	
	\bibitem{3}
	\newblock J. M. Bony,
	\newblock \emph{Equivalence des Diverses Notions de Spectre Singulier Analytique,}
	\newblock S\`{e}minaire Goulaouic-Schwartz, 1976/77, no.3.
	
	\bibitem{4}
	\newblock J. Bros and D. Iagolnitzer,
	\newblock \emph{Support Essentiel et Structure Analytique Des Distributions,}
	\newblock S\`{e}minaire Goulaouic-Lions-Schwartz, 1975/76, no. 18.
	
	\bibitem{5}
	C. R. Crawford, K. F. King, C. J. Ritchie, and J. D. Godwin.
	\textit{Respiratory compensation in projection imaging using a magnification and displacement model.} 
	IEEE Transactions on Medical Imaging, 15 (1996), pp. 327–332.
	
	\bibitem{6}
	L. Desbat, S. Roux, and P. Grangeat.
	\textit{Compensation of some time dependent deformations in tomography.}
	IEEE Transactions on Medical Imaging, 26 (2007), pp. 261–269.
		
	\bibitem{7}
	B. Frigyik, P. Stefanov, and G. Uhlmann. 
	\textit{The X-Ray Transform for a Generic Family of Curves and Weights.} 
	J. Geom. Anal., 18(1):89-108, 2008.
	
	\bibitem{8}
	V. Guillemin, and S. Sternberg. 
	\textit{Some problems in integral geometry and some related problems in microlocal analysis.} 
	Amer. J. Math, 101:915–955, 1979.
	
	\bibitem{9}
	V. Guillemin. 
	\textit{On some results of Gel’fand in integral geometry, in Pseudodifferential operators and applications.}
	Amer. Math. Soc., Providence, RI, 1985.
	
	\bibitem{10}
	V. Guillemin, and S. Sternberg. 
	\textit{Geometric Asymptotics.} American Mathematical Soc., 1990.
	
	\bibitem{11}
	B. N. Hahn, and E. T. Quinto. 
	\textit{Detectable singularities from dynamic Radon data.} 
	SIAM J. Imaging Sciences, 9(3)(2016), pp. 1195–1225.
	
	\bibitem{12}
	B. Hahn.
	\textit{Reconstruction of dynamic objects with affine deformations in dynamic computerized tomography.} 
	J. Inverse Ill-Posed Probl., 22 (2014), pp. 323–339.
	
	\bibitem{13}
	B. N. Hahn.
	\textit{Efficient algorithms for linear dynamic inverse problems with known motion.}
	Inverse	Problems, 30 (2014), pp. 035008, 20.
	
	\bibitem{14}
	L. H\"{o}rmander. 
	\textit{The analysis of linear partial differential operators. III, volume 274. Pseudodifferential operators.} 
	Springer-Verlag, Berlin, 1985. 
	
	\bibitem{15}
	L. H\"{o}rmander. 
	\textit{Fourier Integral Operators, I.} 
	Acta Mathematica, 127 (1971), pp. 79–183. 	
	
	\bibitem{16} 
	\newblock L. H\"ormander,
	\newblock \emph{Uniqueness theorems and wave front sets for solutions of linear differential equations with analytic coefficients,}
	\newblock \emph{Comm. Pure Appl. Math.} \textbf{24} (1971), 671--704.
	
	\bibitem{17}
	A. Homan, and H. Zhou. 
	\textit{Injectivity and stability for a generic class of generalized Radon transforms.} 
	J. Geom. Anal. 27 (2017), no. 2, 1515–1529. 
			
	\bibitem{18}
	A. Katsevich. 
	\textit{Local tomography for the limited-angle problem.}
	J. Math. Anal. Appl., 213 (1997), pp. 160-182.	
	
	\bibitem{19}
	A. Katsevich. 
	\textit{Improved Cone Beam Local Tomography.}
	Inverse Problems, 22 (2006), pp. 627–643.	
		
	\bibitem{20}
	A. Katsevich. 
	\textit{Motion compensated local tomography.}
	Inverse Problems, 24 (2008), 045012.	
	
	\bibitem{21}
	A. Katsevich. 
	\textit{An accurate approximate algorithm for motion compensation in two-dimensional tomography.}
	Inverse Problems, 26 (2010), 065007.
	
	\bibitem{22}
	A. Katsevich, M. Silver, and A. Zamyatin. 
	\textit{Local tomography and the motion estimation problems.}
	SIAM J. Imaging Sci., 4 (2011), pp. 200–219.
		
	\bibitem{23}
	V. P. Krishnan, and E. T. Quinto.
	\textit{Microlocal Analysis in Tomography.}
	In Handbook of Mathematical Methods in Imaging, ed. 2, O. Scherzer, ed., Springer Verlag, 2015. 
	 	 
	\bibitem{24}
	V. Krishnan.
	\textit{A support theorem for the geodesic ray transform on functions.}
	J. Fourier Anal. Appl, 15:515–520, 2009.
	
	\bibitem{25}
	F. Natterer.
	\textit{The mathematics of computerized tomographys.}
	B. G. Teubner, Stuttgart, 1986.
	
	\bibitem{26}
	S. Roux, L. Desbat, A. Koenig, and P. Grangeat.
	\textit{Exact reconstruction in 2d dynamic ct: compensation of time-dependent affine deformations.}
	Physics in Medicine and Biology, 49 (2004), pp. 2169–2182.
	
	\bibitem{27}
	\newblock M. Sato,
	\newblock \emph{Hyperfunctions and Partial Differential Equations,}
	\newblock Proc. Int. Conf. Funct. Anal. Tokyo 1969, 91--4.
	
	\bibitem{28} 
	\newblock J. Sj\"ostrand,
	\newblock Singularit\'es analytiques microlocales,
	\newblock In \emph{ Ast\'erisque, 95, Soc. Math. France,} Paris, volume 95 of Ast\'erisque, (1982), 1--166.
	
	\bibitem{29}
	P. Stefanov, and G. Uhlmann. 
	\textit{Stability estimates for the X-ray transform of tensor fields and boundary rigidity.} 
	Duke Math. J. 123(2004), 445–467.
	 
	\bibitem{30}
	 M. Taylor.  
	\textit{Pseudodifferential Operators.} 
	Princeton University Press, 1981. 
	
	
\end{thebibliography}
\end{document}